\documentclass[11pt]{amsart}
\usepackage[T1]{fontenc}
\usepackage{lmodern,amssymb,amsmath,mathtools}
\usepackage[margin=1.04in]{geometry}
\usepackage{microtype}
\usepackage{xcolor}
\usepackage[colorlinks=true,linkcolor=blue!55!black,citecolor=blue!55!black,urlcolor=blue!55!black]{hyperref}
\hypersetup{pdftitle={An integral test for the maximal operator on variable Lebesgue spaces},pdfauthor={Andrei K. Lerner}}
\newtheorem{theorem}{Theorem}[section]
\newtheorem{lemma}[theorem]{Lemma}
\newtheorem{proposition}[theorem]{Proposition}

\theoremstyle{definition}

\numberwithin{equation}{section}
\newcommand{\R}{\mathbb R}
\newcommand{\av}[1]{\langle #1\rangle}
\newcommand{\Ainf}{\mathcal A_\infty}

\newcommand{\one}{\chi}
\newcommand{\norm}[2]{\lVert #1\rVert_{#2(\cdot)}}
\newcommand{\essinf}{\operatorname*{ess\,inf}}
\newcommand{\esssup}{\operatorname*{ess\,sup}}

\begin{document}
\title[An integral test for the maximal operator]{An integral test for the maximal operator on variable Lebesgue spaces}
\author[A.K. Lerner]{Andrei K. Lerner}
\address[A.K. Lerner]{Department of Mathematics,
Bar-Ilan University, 5290002 Ramat Gan, Israel}
\email{lernera@math.biu.ac.il}

\thanks{The author was supported by Israel Science Foundation grant no. 807/26.}

\begin{abstract}
Let $\mathcal P$ denote the class of the exponents $p(\cdot)$ with $1<p_-\le p_+<\infty$ for which the Hardy--Littlewood maximal operator is bounded on $L^{p(\cdot)}(\mathbb R^n)$.
We obtain a new characterization of $\mathcal P$ in terms of a single integral condition on the exponent itself. The criterion requires neither a separate Muckenhoupt condition nor testing on families of pairwise disjoint cubes. As an application, we answer an open question by showing that if $p(\cdot)\in {\mathcal P}$, then $p(\cdot)+c\in {\mathcal P}$ for every constant $c>0$.
\end{abstract}

\keywords{Maximal operator, variable Lebesgue spaces,
integral criteria for boundedness, perturbations of variable exponents.}
\subjclass[2020]{42B25, 46E30}

\maketitle

\section{Introduction}
Let $p(\cdot):\R^n\to[1,\infty)$ be a measurable function. The space $L^{p(\cdot)}(\R^n)$ consists of measurable functions $f$ with finite Luxemburg norm
\[
 \norm f p=\inf\left\{\lambda>0:
 \int_{\R^n}\left(\frac{|f(x)|}{\lambda}\right)^{p(x)}dx\le1\right\}.
\]
For a cube $Q$ and a locally integrable function $f$, write $\av f_Q=|Q|^{-1}\int_Q f$. The Hardy--Littlewood maximal operator is defined by
\[
 Mf(x)=\sup_{Q\ni x}\av{|f|}_Q,
\]
where the supremum is taken over all cubes containing $x$ with sides parallel to the coordinate axes.

Our aim is to characterize the exponents for which $M$ is bounded on $L^{p(\cdot)}(\R^n)$. For sufficient conditions in terms of log-H\"older continuity and for the history of the problem, we refer to the monographs \cite{CF,DHHR}. Write $p_-=\essinf p(\cdot)$ and $p_+=\esssup p(\cdot)$. Throughout the paper we consider the range $1<p_-\le p_+<\infty$.

Diening \cite{Diening} proved that boundedness of $M$ is equivalent to uniform boundedness of the operators
\[
 f\longmapsto \sum_{Q\in\mathcal F}\av f_Q\,\one_Q
\]
over finite families $\mathcal F$ of pairwise disjoint cubes. Shalukhina \cite{Shalukhina} recently gave a shorter proof of this result. Testing on a single cube gives the variable Muckenhoupt condition $A_{p(\cdot)}$:
\[
 \sup_Q\frac{\norm{\one_Q}p\,\norm{\one_Q}{p'}}{|Q|}<\infty,
 \qquad p'(x)=\frac{p(x)}{p(x)-1}.
\]
This condition is necessary but not sufficient \cite[Example~4.51]{CF}. It does, however, control the boundedness of $M$ on compact sets \cite{Kopaliani,Shifts}.
This suggests looking for a characterization of the form $A_{p(\cdot)}\cap B$, where $B$ describes the behavior of $p(\cdot)$ at infinity. For example, $A_{p(\cdot)}$ together with Nekvinda's integral decay condition \cite{Nekvinda} is sufficient for boundedness \cite[Theorem~4.52]{CF}. This decay condition controls the deviation of $p(\cdot)$ from a fixed constant at infinity, but it is not necessary. Replacing it by the weaker condition $\mathcal U_\infty$, one obtains a necessary and sufficient condition: $M$ is bounded if and only if $p(\cdot)\in A_{p(\cdot)}\cap\mathcal U_\infty$ \cite{UnionOld}.

A different characterization was obtained in \cite{Ainf}: $M$ is bounded if and only if both $p(\cdot)$ and $p'(\cdot)$ satisfy the variable analogue of the $A_\infty$ condition. A recent characterization in \cite{GOR} applies also to weighted variable Lebesgue spaces. In the unweighted case it requires $A_{p(\cdot)}$ and the existence of a positive integrable function $h$ such that $h^{1/p(\cdot)}$ and $h^{1/p'(\cdot)}$ belong to $A_1$.

The four characterizations \cite{Diening,UnionOld,Ainf,GOR} described above express boundedness in quite different terms. One reason to seek another
characterization is to understand how boundedness behaves under perturbations of the exponent. Consider the following two questions.

\vskip 2mm
\noindent\textbf{Question 1.} \emph{Suppose that $M$ is bounded on both $L^{p_1(\cdot)}$ and $L^{p_2(\cdot)}$. Is $M$ bounded on $L^{p_1(\cdot)+p_2(\cdot)}$?}

\vskip 2mm
\smallskip
\noindent\textbf{Question 2.} \emph{Let $c>0$. If $M$ is bounded on $L^{p(\cdot)}$, is it also bounded on $L^{p(\cdot)+c}$?}

\vskip 2mm
Question~2 is a special case of Question~1 when $c>1$; for arbitrary $c>0$ it is posed explicitly in the monograph of Cruz-Uribe and Fiorenza \cite[Problem~A.10; see also Question~4.75]{CF}.
To our knowledge, neither question has previously been answered in full. The answers are not apparent from the characterizations cited above.

Our main result characterizes boundedness by a single integral condition on the exponent $p(\cdot)$ itself. It requires neither a separate
$A_{p(\cdot)}$ condition nor testing on families of pairwise disjoint cubes. Using this criterion, we answer both questions affirmatively.

To state the result, for a cube $Q$ and $s>0$ define
\[
 F_Q(s)=\max_{r\in\{p,p'\}}\frac{1}{|Q|}\int_Q s^{\,\frac{r(y)}{r_Q}-1}\,dy,
 \qquad r_Q=\av r_Q .
\]
By Jensen's inequality, $F_Q(s)\ge1$ for every $s>0$.

\begin{theorem}\label{main}
Let $1<p_-\le p_+<\infty$, and let $Q_v(x)$ be the cube with center $x$ and side length $v$. The maximal operator is bounded on $L^{p(\cdot)}(\R^n)$ if and only if, for some $C>1$,
\[
 \int_{\R^n}\int_0^\infty
 \left|\left\{0<s<v^{-n}: F_{Q_v(x)}(s)>C\right\}\right|
 \,\frac{dv}{v}\,dx<\infty .
\]
The inner bars denote one-dimensional Lebesgue measure.
\end{theorem}

If $p(\cdot)$ is constant on $Q$, then $F_Q\equiv1$ and $Q$ contributes nothing. Otherwise $F_Q(s)$ increases from $1$ to infinity as $s$ moves away from $1$ toward either $0$ or $\infty$.
Within $(0,1/|Q|)$, the set where $F_Q(s)>C$ therefore consists of an interval next to zero and possibly a second interval ending at
$1/|Q|$ (see Lemma~\ref{level-intervals} for details). Their total length measures the range of $s$ on which the variation of the exponent on $Q$ makes at least one of the two power averages exceed $C$. Finiteness thus limits the cumulative effect of these oscillations over all locations and scales. For example, the same nonconstant pattern cannot recur on infinitely many disjoint translates of a cube: the subcubes of each copy contribute at least a fixed positive amount.

To apply Theorem~\ref{main}, it suffices to estimate $F_Q(s)$ from above. In particular, one can transfer boundedness from one exponent to another
by comparing their functions $F_Q(s)$ on each cube. To illustrate this, write $F_Q(s;p)$ to indicate the dependence on the exponent. By
Proposition~\ref{shift}, for every $c,s>0$ and $Q$ we have
\[
 F_Q(s;p+c)\le F_Q(s;p).
\]
Thus, for the same $C$, the integral in Theorem~\ref{main} cannot increase when $p(\cdot)$ is replaced by $p(\cdot)+c$. This gives an affirmative answer to Question~2.

The paper is organized as follows. Section~2 collects the preliminary results, and Section~3 proves Theorem~\ref{main}. In
Section~\ref{examples} we answer Questions~1 and~2. We also determine, in terms of the heights and widths of disjoint smooth bumps, exactly
when adding them to an exponent preserves boundedness of $M$.

\section{Preliminaries}
\subsection{Norms and two boundedness results}
Let $1<r_-\le r_+<\infty$ and write
\[
 \rho_r(f)=\int_{\R^n}|f(x)|^{r(x)}\,dx .
\]
The unit-ball property of the Luxemburg norm \cite[Lemma~3.2.4]{DHHR}
states that
\[
 \norm f r\le1 \quad\Longleftrightarrow\quad \rho_r(f)\le1 .
\]
For $f\in L^{r(\cdot)}$, write $N=\norm f r$. By
\cite[Lemma~3.2.5]{DHHR},
\[
 \min\{N^{r_-},N^{r_+}\}\le\rho_r(f)\le\max\{N^{r_-},N^{r_+}\}.
\]
In particular, a uniform bound on the modular implies a uniform bound
on the norm. If $f\ne0$, then $\norm{f/N}r=1$, so both bounds for
$\rho_r(f/N)$ are equal to one. Hence $\rho_r(f/N)=1$.

The variable H\"older inequality \cite[Lemma~3.2.20]{DHHR} states that
\[
 \int_{\R^n}|f(x)g(x)|\,dx\le 2\norm f r\,\norm g {r'} .
\]
The associate space of $L^{r(\cdot)}$ is $L^{r'(\cdot)}$, with equivalent norms. More precisely, by \cite[Theorem~3.2.13 and Corollary~3.2.14]{DHHR},
\[
 \tfrac12\norm g {r'}
 \le \sup_{\norm f r\le1}\int_{\R^n}|f(x)g(x)|\,dx
 \le 2\norm g {r'} .
\]
Taking $r=p$ and $g=\one_Q$ in this formula shows why boundedness of $M$ implies
\[
 \sup_Q A(Q)<\infty,\qquad
 A(Q)=\frac{\norm{\one_Q}p\,\norm{\one_Q}{p'}}{|Q|}.
\]
Indeed, $\av{|f|}_Q\one_Q\le Mf$. Taking norms and then the supremum over $\norm f p\le1$ yields $A(Q)/2\le\|M\|_{L^{p(\cdot)}\to L^{p(\cdot)}}$ for every $Q$.

Several results will be used in the proof of Theorem \ref{main}. First, boundedness of $M$ on $L^{r(\cdot)}$ is equivalent to boundedness on $L^{r'(\cdot)}$; see \cite[Theorem~8.1]{Diening} and also
\cite[Corollary~1.5]{Dual} for a different proof. In addition, boundedness on $L^{r(\cdot)}$ gives boundedness on $L^{r(\cdot)/s}$ for $s>1$ sufficiently close to one (see, e.g., \cite[Theorem~1.2]{LO}).
Second, we use the following characterization from \cite[Theorem~1.3]{Ainf}. We give its formulation to specify exactly what will be needed in the sufficiency argument.

\begin{theorem}[Variable $\Ainf$ criterion]\label{blackbox}
For $1<p_-\le p_+<\infty$, boundedness of $M$ on $L^{p(\cdot)}$ is equivalent to the following property for each $r(\cdot)=p(\cdot),p'(\cdot)$: there exist $0<\lambda<1$ and $B<\infty$ such that
\[
 \norm{\sum_Q u_Q\one_Q}r
 \le B\norm{\sum_Q u_Q\one_{E_Q}}r
\]
for every finite family of pairwise disjoint cubes, all scalars $u_Q\ge0$, and all measurable $E_Q\subset Q$ with $|E_Q|\ge\lambda|Q|$.
\end{theorem}
The finite-family version gives the countable version by monotone convergence. We will verify the property with $\lambda=1/2$.

\subsection{From arithmetic means to harmonic means}
The statement of Theorem~\ref{main} uses the arithmetic mean $r_Q$.
In the proof it is convenient to use the harmonic mean. For an
exponent $r(\cdot)$, write
\[
 h_{r,Q}=\left\langle\frac1r\right\rangle_Q^{-1},
 \qquad
 d_{r,Q}(y)=\frac{r(y)}{h_{r,Q}}-1.
\]
In place of $F_Q$ we shall use the function
\[
 \Phi_Q(t)=
 \left\langle\max\big\{t^{d_{p,Q}},t^{d_{p',Q}}\big\}\right\rangle_Q,
 \qquad t>0.
\]

Let us explain the choice of the harmonic mean and why the maximum
is taken before averaging. Since $1/h_{p,Q}+1/h_{p',Q}=1$, we have
\begin{equation}\label{opposite}
 \frac{d_{p,Q}(y)}{p(y)}
 =-\frac{d_{p',Q}(y)}{p'(y)}
 =\frac{1}{h_{p,Q}}-\frac{1}{p(y)}.
\end{equation}
For $t>1$, the larger of the two powers is $t^{d_{p,Q}(y)}$ where
$1/p(y)<1/h_{p,Q}$, and $t^{d_{p',Q}(y)}$ where
$1/p(y)>1/h_{p,Q}$. For $0<t<1$ the roles are reversed. Thus the
larger power need not be the same throughout $Q$. Taking the maximum
at each point allows us to estimate the absolute deviation of $1/p(y)$
from its average on $Q$, whichever sign the deviation has. Indeed,
by \eqref{opposite} and the fact that $p(y),p'(y)>1$, we have
\[
 \max\big\{t^{d_{p,Q}(y)},t^{d_{p',Q}(y)}\big\}
 \ge
 \exp\!\left(
   |\log t|\left|\frac1{p(y)}-\frac1{h_{p,Q}}\right|
 \right).
\]
Averaging over $Q$ and applying Jensen's inequality to the exponential,
we obtain
\begin{equation}\label{reciprocal-osc}
 |\log t|
 \left\langle\left|\frac1p-\frac1{h_{p,Q}}\right|\right\rangle_Q
 \le \log\Phi_Q(t).
\end{equation}
We shall use \eqref{reciprocal-osc} to compare the average of
$1/p(\cdot)$ on a cube with its average on a subset of comparable
measure.

\begin{lemma}\label{centers}
For every $C>1$ there is $D>1$, depending only on $C,p_-,p_+$, such that for every cube $Q$ and every $t>0$,
\[
 F_Q(t)\le C\ \Longrightarrow\ \Phi_Q(t)\le D,
 \qquad
 \Phi_Q(t)\le C\ \Longrightarrow\ F_Q(t)\le D.
\]
\end{lemma}
\begin{proof}
Fix $r=p$ or $p'$. Put $f=(r/r_Q-1)\log t$ and $g=d_{r,Q}\log t$. The difference between the two normalizations can be computed in either direction:
\[
 g-f=-r\langle f/r\rangle_Q
     =\frac{r}{r_Q}\langle g\rangle_Q .
\]
Indeed, $\langle f/r\rangle_Q=(1/r_Q-1/h_{r,Q})\log t$, whereas
$\langle g\rangle_Q=(r_Q/h_{r,Q}-1)\log t$.

Suppose first that $F_Q(t)\le C$. Then $\langle f\rangle_Q=0$ and $\langle e^f\rangle_Q\le C$. For any real function $u$ with $\langle u\rangle_Q=0$, its positive part $u^+=\max\{u,0\}$ has the same integral as its negative part. Jensen's inequality and $e^{u^+}\le 1+e^u$ therefore give
\begin{equation}\label{zero}
 \langle |u|\rangle_Q=2\langle u^+\rangle_Q\le2\log(1+C).
\end{equation}
Apply this to $f$. The first expression for $g-f$ gives $|g-f|\le (r_+/r_-)\langle|f|\rangle_Q$, so $\langle e^g\rangle_Q$ is bounded for both choices of $r$. Their sum bounds $\Phi_Q(t)$.

Conversely, suppose $\Phi_Q(t)\le C$. By \eqref{opposite} and \eqref{reciprocal-osc},
$\langle|g|\rangle_Q\le r_+\log C$. The second expression for $g-f$ now gives a uniform bound for $|g-f|$. Since $\langle e^g\rangle_Q\le\Phi_Q(t)\le C$, we obtain a bound for $\langle e^f\rangle_Q$ for each $r$, and hence for $F_Q(t)$.
\end{proof}

By Lemma~\ref{centers}, the integral in Theorem~\ref{main} is finite
for some $C>1$ if and only if the same is true with $\Phi_Q$ in place
of $F_Q$.

The following lemma describes the sets whose lengths occur in the
integral criterion. We shall also use it in the proof of
Lemma~\ref{small}: to bound $\Phi_Q$ on an interval containing $1$,
it suffices to bound it at the two endpoints.

\begin{lemma}\label{level-intervals}
For every cube $Q$, the functions $F_Q$ and $\Phi_Q$ are continuous
on $(0,\infty)$, equal one at $t=1$, and are nonincreasing on $(0,1]$
and nondecreasing on $[1,\infty)$. Consequently, for every $C>1$,
each of the sets
\[
 \{0<t<1/|Q|:F_Q(t)>C\}
 \qquad\text{and}\qquad
 \{0<t<1/|Q|:\Phi_Q(t)>C\}
\]
is a union of at most two intervals.
\end{lemma}

\begin{proof}
Fix a cube $Q$. The functions $u\mapsto F_Q(e^u)$ and
$u\mapsto\Phi_Q(e^u)$ are convex on $\mathbb R$. Indeed, in both
definitions we take averages and maxima of exponentials of linear
functions of $u$, and both operations preserve convexity.

Both functions take the value one at $u=0$, and this is their minimum:
for $F_Q$ this follows from Jensen's inequality, and for $\Phi_Q$
from \eqref{reciprocal-osc}. A finite convex function with a minimum
at zero is continuous, nonincreasing to the left of zero, and
nondecreasing to its right. Since $\log t$ is increasing and vanishes
at $t=1$, the asserted continuity and monotonicity follow.

If $p(\cdot)$ is nonconstant on $Q$, the minimum is strict in both
cases, by strict Jensen's inequality for the term $r=p$ in $F_Q$ and
by \eqref{reciprocal-osc} for $\Phi_Q$. A convex function with a
strict minimum at zero tends to infinity as $u\to\pm\infty$.
Thus, for each of $F_Q$ and $\Phi_Q$, the set of $t>0$ where the
function is at most $C$ is a closed interval containing $1$.
Its complement intersects $(0,1/|Q|)$ in at most two intervals.
If $p(\cdot)$ is constant on $Q$, both functions are identically
one and the two sets in the statement are empty.
\end{proof}

In the proof of Theorem~\ref{main}, we use the following notation:
\[
 \Lambda_C(Q)=|\{0<t<1/|Q|:\Phi_Q(t)>C\}|,\qquad
 I_C=\int_{\R^n}\int_0^\infty\Lambda_C(Q_v(x))\,\frac{dv}{v}\,dx .
\]
To check measurability when the cube varies, observe that $(x,v,t)\mapsto\Phi_{Q_v(x)}(t)$ is continuous for $v,t>0$. Indeed, along a convergent sequence the cube indicators converge almost everywhere and are supported in a common bounded set. Dominated convergence applies first to the averages of $1/p$ and then to the defining integral for $\Phi$. Consequently $\Lambda_C(Q_v(x))$ is the integral in $t$ of the indicator of the open set $\{(x,v,t):v>0,\ 0<t<v^{-n},\ \Phi_{Q_v(x)}(t)>C\}$. Fubini's theorem for nonnegative functions gives measurability in $(x,v)$. The same argument, using the averages of $p$ and $p'$, applies to the integral with $F_Q$.

\subsection{The endpoint and small cubes}
The value $\Phi_Q(1/|Q|)$ has a familiar meaning: it bounds the modulars of two normalized characteristic functions. We record the relation with the Muckenhoupt condition in a form that also applies to an individual cube. For the harmonic-mean norm estimates see \cite[Theorem~4.5.7]{DHHR} or \cite[Proposition~4.66]{CF}.

\begin{lemma}\label{endpoint}
On each cube $Q$, a bound for $A(Q)$ gives a bound for $\Phi_Q(1/|Q|)$ and conversely. The bounds depend only on one another and on $p_-,p_+$. The same conclusion holds for $F_Q(1/|Q|)$.
\end{lemma}

\begin{proof}
Write $V=|Q|$. We first relate $\Phi_Q(1/V)$ to the modulars of two
normalized characteristic functions. For $r=p$ or $r=p'$, the identity
$d_{r,Q}(x)=r(x)/h_{r,Q}-1$ gives
\[
 \rho_r\bigl(V^{-1/h_{r,Q}}\one_Q\bigr)
 =\int_Q V^{-r(x)/h_{r,Q}}\,dx
 =\frac1V\int_Q (1/V)^{d_{r,Q}(x)}\,dx.
\]
Since $\max\{a,b\}\le a+b\le 2\max\{a,b\}$ for $a,b\ge0$, it follows
from the definition of $\Phi_Q$ that
\[
 \Phi_Q(1/V)
 \le \rho_p\bigl(V^{-1/h_{p,Q}}\one_Q\bigr)
    +\rho_{p'}\bigl(V^{-1/h_{p',Q}}\one_Q\bigr)
 \le 2\Phi_Q(1/V).
\]
Thus a bound for $\Phi_Q(1/V)$ is equivalent to a bound for both
modulars. By the norm--modular inequalities in Section~2.1, this is
equivalent to a bound for both ratios
\[
 \frac{\norm{\one_Q}p}{V^{1/h_{p,Q}}}
 \qquad\text{and}\qquad
 \frac{\norm{\one_Q}{p'}}{V^{1/h_{p',Q}}}.
\]

We now compare these ratios with $A(Q)$. Since
$1/h_{p,Q}+1/h_{p',Q}=1$, their product is
\[
 \frac{\norm{\one_Q}p}{V^{1/h_{p,Q}}}
 \frac{\norm{\one_Q}{p'}}{V^{1/h_{p',Q}}}
 =\frac{\norm{\one_Q}p\,\norm{\one_Q}{p'}}{V}
 =A(Q).
\]
Hence a bound for both ratios gives a bound for $A(Q)$.

Now suppose that $A(Q)$ is bounded. We want to bound each ratio from
above. Since their product is $A(Q)$, each ratio equals $A(Q)$ divided
by the other. We therefore show that both ratios are at least $1/2$;
it will follow that each is at most $2A(Q)$.

Fix $r=p$ or $r=p'$. By Jensen's inequality,
\[
 V^{1/h_{r,Q}}
 =\exp\!\left(\frac1V\int_Q\log V^{1/r(x)}\,dx\right)
 \le \frac1V\int_Q V^{1/r(x)}\,dx
 =\int_Q V^{-1/r'(x)}\,dx.
\]
The function $V^{-1/r'(\cdot)}\one_Q$ has $L^{r'(\cdot)}$ norm one,
since its modular is $\int_Q V^{-1}\,dx=1$. Therefore the variable
H\"older inequality gives
\[
 \int_Q V^{-1/r'(x)}\,dx\le 2\norm{\one_Q}r.
\]
This proves that each ratio is at least $1/2$. Since their product
is $A(Q)$, each is at most $2A(Q)$. Thus a bound for $A(Q)$ also
gives a bound for $\Phi_Q(1/V)$.

Finally, Lemma~\ref{centers}, applied at $t=1/V$, gives the same
equivalence with $F_Q(1/V)$ in place of $\Phi_Q(1/V)$.
\end{proof}

\begin{lemma}\label{small}
There exists $0<\eta\le1$, depending only on $p_-,p_+$, such that on every cube $Q$ with $|Q|<1$,
\[
 \Phi_Q(t)\le\Phi_Q(1/|Q|)\qquad
 \bigl(|Q|^\eta\le t\le 1/|Q|\bigr).
\]
\end{lemma}
\begin{proof}
Write $V=|Q|$. By monotonicity on either side of $1$, it suffices to prove the inequality at $t=V^\eta$. Choose $\eta$ so that $\eta(r(y)-1)\le1$ for $r=p,p'$ almost everywhere. By \eqref{opposite}, $-d_{r,Q}=(r-1)d_{r',Q}$. Where $d_{r,Q}<0$, the conjugate exponent therefore gives
\[
 V^{\eta d_{r,Q}}=(1/V)^{\eta(r-1)d_{r',Q}}
 \le (1/V)^{d_{r',Q}}.
\]
Where $d_{r,Q}\ge0$, the left-hand side is at most one. The integrand in $\Phi_Q(1/V)$ is at least one at each point and thus bounds both powers at $V^\eta$. Averaging proves the lemma.
\end{proof}

\section{Proof of Theorem~\ref{main}}
By Lemma~\ref{centers} it suffices to prove that $M$ is bounded on $L^{p(\cdot)}$ if and only if $I_C<\infty$ for some $C>1$. We work with $\Phi_Q$ throughout the proof.

\subsection{Necessity}
For fixed $t$ and side length $v$, we first bound the measure of the centers $x$ where $\Phi_{Q_v(x)}(t)$ is large. The key is that the bound grows more slowly than $1/t$ as $t$ tends to zero. The following lemma gives this improvement before we integrate in either variable.

\begin{lemma}\label{packing}
Suppose $1<r_-\le r_+<\infty$ and $M$ is bounded on $L^{r(\cdot)}$. There exist $K>1$ and $0<\kappa<1$ such that, for $0<t<1$ and every finite disjoint family $\mathcal F$ of cubes,
\[
 \langle t^{d_{r,Q}}\rangle_Q>K\quad\text{for all }Q\in\mathcal F
 \qquad\Longrightarrow\qquad
 \sum_{Q\in\mathcal F}|Q|<t^{-\kappa}.
\]
\end{lemma}
\begin{proof}
By duality and the self-improvement in Section~2.1, $M$ is bounded on $L^{q(\cdot)}$, where $q(\cdot)$ is the conjugate of $r'(\cdot)/s$ for some $1<s<(r')_-$. The reason to use this exponent is the identity
\[
 \frac1q=\frac{s}{r}-(s-1),\qquad
 d_{r,Q}(y)=\bigl(1-(1-1/s)r(y)\bigr)d_{q,Q}(y).
\]
Indeed, $\langle1/q\rangle_Q-1/q=s(\langle1/r\rangle_Q-1/r)$, and $d_{r,Q}=r(\langle1/r\rangle_Q-1/r)$. The factor multiplying $d_{q,Q}(y)$ is positive since $s<(r')_-$, and it is at most $\kappa=1-(1-1/s)r_-<1$.

Suppose that $E=\bigcup_{Q\in\mathcal F}Q$ has measure at least $t^{-\kappa}$ and put $u=1/|E|$. Then $u\le t^\kappa\le t^{\,1-(1-1/s)r(y)}\le1$. Raising a number in $[u,1]$ to a real power gives at most the larger of its two endpoint values. Thus
\[
 t^{d_{r,Q}(y)}\le \max\{1,u^{d_{q,Q}(y)}\}.
\]
To average the right-hand side over $E$, take $f=u^{1/q(\cdot)}\one_E$. Its modular is $u|E|=1$, so boundedness of $M$ gives $\rho_q(Mf)\le B$ independently of the family. On $Q$, Jensen's inequality gives 
$$Mf\ge\langle u^{1/q}\rangle_Q\ge u^{1/h_{q,Q}}.$$ 
Since the cubes are disjoint, we obtain
\[
 \sum_{Q\in\mathcal F}\frac{|Q|}{|E|}\,\langle t^{d_{r,Q}}\rangle_Q
 \le 1+\sum_{Q\in\mathcal F}\int_Q u^{q(y)/h_{q,Q}}\,dy
 \le 1+\rho_q(Mf)\le1+B .
\]
The weights on the left sum to one. Hence at least one cube has average at most $1+B$, which proves the lemma with $K=1+B$.
\end{proof}

Now assume that $M$ is bounded on $L^{p(\cdot)}$. Apply Lemma~\ref{packing} to $p(\cdot)$ and to $p'(\cdot)$, using dual boundedness in the latter case. Let $K$ and $\kappa$ be the larger of the two constants and the larger of the two exponents, respectively. If $\Phi_Q(t)>2K$, then at least one of $\langle t^{d_{p,Q}}\rangle_Q$ and $\langle t^{d_{p',Q}}\rangle_Q$ exceeds $K$, because their sum bounds $\Phi_Q(t)$. Splitting any finite disjoint family according to which average exceeds $K$ gives
\[
 \sum_{Q\in\mathcal F}|Q|\le 2t^{-\kappa}
 \qquad\text{if }0<t<1\text{ and }\Phi_Q(t)>2K\text{ on every }Q\in\mathcal F .
\]

Choose $C\ge2K$ also large enough that $\Phi_Q(1/|Q|)\le C$ on every cube. Such a choice is possible by Section~2.1 and Lemma~\ref{endpoint}. For $v,t>0$ define
\[
 E_{v,t}=\{x\in\R^n:\Phi_{Q_v(x)}(t)>C\}.
\]
For fixed $z\in[0,v)^n$, the cubes $Q_v(z+vk)$, $k\in\mathbb Z^n$, have disjoint interiors. Each has volume $v^n$. Applying the last estimate to finite subfamilies and then passing to their union shows that, for $0<t<1$,
\[
 v^n\sum_{k\in\mathbb Z^n}\one_{E_{v,t}}(z+vk)\le 2t^{-\kappa}.
\]
Integrating over $z\in[0,v)^n$ counts every center exactly once, since the translates of $[0,v)^n$ partition~$\R^n$. We obtain
\begin{equation}\label{slice}
 |E_{v,t}|=\int_{[0,v)^n}\sum_{k\in\mathbb Z^n}
       \one_{E_{v,t}}(z+vk)\,dz\le 2t^{-\kappa}.
\end{equation}

It remains to identify which side lengths contribute. In $I_C$ the parameter satisfies $0<t<v^{-n}$. If $v<1$, Lemma~\ref{small} gives $\Phi_{Q_v(x)}(t)\le C$ whenever $v^{n\eta}\le t\le v^{-n}$. Thus a nonempty $E_{v,t}$ requires $t<v^{n\eta}<1$, or $v>t^{1/(n\eta)}$. If $v\ge1$, then $t<v^{-n}\le1$ already implies $t^{1/(n\eta)}<1\le v$. In both cases it suffices to integrate over $0<t<1$ and $t^{1/(n\eta)}<v<t^{-1/n}$. By Fubini's theorem and \eqref{slice},
\[
 \begin{aligned}
 I_C&=\int_0^\infty\int_0^{t^{-1/n}}|E_{v,t}|\,\frac{dv}{v}\,dt\\
 &\le 2\int_0^1 t^{-\kappa}\int_{t^{1/(n\eta)}}^{t^{-1/n}}
       \frac{dv}{v}\,dt
 =\frac{2(1+\eta^{-1})}{n}\int_0^1 t^{-\kappa}\log(1/t)\,dt<\infty .
 \end{aligned}
\]
This proves necessity.

\subsection{Sufficiency}
Suppose $I_{C_0}<\infty$ for some $C_0>1$. We will use the fact that $I_C\to0$ as $C\to\infty$. Indeed, for each cube $Q$ the sets defining $\Lambda_C(Q)$ decrease to the empty set and have measure at most $1/|Q|$. Hence $\Lambda_C(Q)\to0$, and dominated convergence applies since $I_{C_0}<\infty$.

\emph{Step 1. An estimate for disjoint cubes.}
We first show that there is a level $K>1$ such that for every finite disjoint family $\mathcal F$,
\begin{equation}\label{budget}
 \sum_{Q\in\mathcal F}|Q|\Lambda_K(Q)\le\tfrac14 .
\end{equation}
To relate the length on one cube to the integral over centers and side lengths, suppose $Q\subset R$ and $|R|\le k|Q|$, where $k=3^n$. Write $a=1/p$ and $\beta=\max\{p_+,(p')_+\}$. By \eqref{reciprocal-osc},
\[
 |a_Q-a_R|\,|\log t|
 \le k\,\langle|a-a_R|\rangle_R\,|\log t|
 \le k\log\Phi_R(t).
\]
For $r=p,p'$ we have $|d_{r,Q}-d_{r,R}|=r|a_Q-a_R|$ by \eqref{opposite}. Therefore $t^{d_{r,Q}}\le\Phi_R(t)^{k\beta}t^{d_{r,R}}$ on $Q$. Taking the maximum over $r$ and averaging first on $Q$, then bounding that average by $k$ times the average on $R$, gives
\[
 \Phi_Q(t)\le k\,\Phi_R(t)^{1+k\beta}.
\]
Also $\Phi_R(t)\le k^\beta\Phi_R(t/k)$ because $|d_{r,R}|\le\beta$. For any $C>1$, take $K=k(k^\beta C)^{1+k\beta}$. If $t$ is counted by $\Lambda_K(Q)$, the last two estimates give $\Phi_R(t/k)>C$. Since $t/k<1/(k|Q|)\le1/|R|$, this proves $\Lambda_K(Q)\le k\Lambda_C(R)$.

For each $x\in Q$ and $2\ell_Q<v<3\ell_Q$, where $\ell_Q$ is the side length of $Q$, the cube $R=Q_v(x)$ contains $Q$ and has volume at most $k|Q|$. Integrating the length estimate and summing in $Q$ thus gives
\[
 \sum_{Q\in\mathcal F}|Q|\Lambda_K(Q)
 \le \frac{k}{\log(3/2)}
       \sum_{Q\in\mathcal F}\int_Q\int_{2\ell_Q}^{3\ell_Q}
         \Lambda_C(Q_v(x))\,\frac{dv}{v}\,dx
 \le \frac{k}{\log(3/2)}I_C .
\]
The last inequality uses only disjointness in the center variable $x$. Choose $C$ so large that the last expression is at most $1/4$, and fix the corresponding $K$. This proves \eqref{budget}.

\emph{Step 2. Passing from a cube to a subset of half its measure.}
Fix $r=p$ or $p'$ and take a family, scalars $u_Q$ and sets $E_Q$ as in Theorem~\ref{blackbox}, with $|E_Q|\ge|Q|/2$. Omit zero terms. By homogeneity normalize the norm of $\sum_Q u_Q\one_Q$ to one. Its modular is then one, so the numbers
\[
 m_Q=\int_Q u_Q^{r(y)}\,dy\qquad\text{satisfy}\qquad \sum_Q m_Q=1 .
\]
We will show that a fixed part of this modular remains on the sets $E_Q$.

For each $Q$, consider the interval $m_Q/(2|Q|)<t<m_Q/|Q|$. It lies in $(0,1/|Q|)$ since $m_Q\le1$. If $\Phi_Q(t)>K$ throughout this interval, its length gives $m_Q\le 2|Q|\Lambda_K(Q)$. By \eqref{budget}, all such cubes together carry at most half the modular. It suffices to treat the other cubes, whose $m_Q$ sum to at least $1/2$.

On each of these cubes choose $t$ in the indicated interval with $\Phi_Q(t)\le K$. The inequality $t|Q|\le m_Q$ allows us to compare $u_Q$ with $t^{1/h_{r,Q}}$. Indeed, since $K>1$ and $r(y)>1$,
\[
 \int_Q (t^{1/h_{r,Q}}/K)^{r(y)}\,dy
 \le K^{-1}t|Q|\Phi_Q(t)\le m_Q=\int_Q u_Q^{r(y)}\,dy .
\]
The integral is increasing in the scalar, so $u_Q\ge t^{1/h_{r,Q}}/K$. On $E_Q$, \eqref{opposite} and \eqref{reciprocal-osc} give
\[
 |\langle d_{r,Q}\rangle_{E_Q}\log t|
 \le 2\beta\langle|a-a_Q|\rangle_Q\,|\log t|\le2\beta\log K .
\]
Jensen's inequality therefore gives $\langle t^{d_{r,Q}}\rangle_{E_Q}\ge K^{-2\beta}$. Together with $|E_Q|\ge|Q|/2$ and $t>m_Q/(2|Q|)$, this yields
\[
 \int_{E_Q}u_Q^{r(y)}\,dy
 \ge K^{-\beta}t|E_Q|\langle t^{d_{r,Q}}\rangle_{E_Q}
 \ge \tfrac14 K^{-3\beta}m_Q .
\]
Summing over these cubes implies $\rho_r(\sum_Q u_Q\one_{E_Q})\ge K^{-3\beta}/8$. By Section~2.1 the norm of this sum is bounded below by a positive constant, while the norm on the full cubes is one. This proves the norm inequality of Theorem~\ref{blackbox} for both exponents and completes the proof of Theorem~\ref{main}. 

\section{Two applications}\label{examples}

\subsection{Addition of exponents}
We prove a single proposition that answers Questions~1 and~2 from the introduction.
\begin{proposition}\label{addition}
Suppose that $1<(p_i)_-\le(p_i)_+<\infty$ and that $M$ is bounded on $L^{p_i(\cdot)}$ for $i=1,2$. Let $t>0$ and put $P(\cdot)=p_1(\cdot)+t p_2(\cdot)$. If $P(\cdot)\ge p_2(\cdot)$ almost everywhere, then $M$ is bounded on $L^{P(\cdot)}$.
\end{proposition}

For Question~1 take $t=1$. The hypothesis $p_1(\cdot)+p_2(\cdot)\ge p_2(\cdot)$ then holds automatically.

For Question~2 take $p_1(\cdot)=p(\cdot)$, $p_2(\cdot)\equiv p_-$, and $t=c/p_-$. The constant exponent on which we use boundedness is $p_->1$, while the amount added to $p(\cdot)$ is $t p_-=c$. Since $p(\cdot)+c\ge p_-$ almost everywhere, the proposition applies for every $c>0$.

\begin{proof}
We use the harmonic means from Section~2. Fix a cube $Q$ and $s>0$. Write $h_r=\langle 1/r\rangle_Q^{-1}$. For any $b>1$, put
\[
 W(r,b)=\max\{s^{r/b-1},s^{r'/b'-1}\},\qquad b'=b/(b-1).
\]
Thus $\langle W(r,h_r)\rangle_Q=\Phi_Q(s;r)$, where the last argument specifies the exponent used in $\Phi_Q$. The two powers have exponents of opposite signs, since
\[
 \frac{r'}{b'}-1=-\frac{r/b-1}{r-1}.
\]
In particular, $W(r,b)\ge1$.

We first compare the two powers for $P$ with those for $p_1$ and $p_2$. Choose $b=h_{p_1}+t h_{p_2}$, $\alpha_1=h_{p_1}/b$ and $\alpha_2=t h_{p_2}/b$. Since $\alpha_1+\alpha_2=1$, we obtain
\[
 \frac{P}{b}-1=\sum_{i=1}^2\alpha_i
       \left(\frac{p_i}{h_{p_i}}-1\right).
\]
Multiplying this equality by $-1/(P-1)$ gives the corresponding identity for the conjugates:
\[
 \frac{P'}{b'}-1=\sum_{i=1}^2\alpha_i\,
       \frac{p_i-1}{P-1}\left(\frac{p_i'}{h_{p_i}'}-1\right).
\]
Here $0<(p_i-1)/(P-1)\le1$ by assumption. Taking a power between zero and one of $s^{p_i'/h_{p_i}'-1}$ leaves it bounded by $W(p_i,h_{p_i})$, since this maximum is at least one. The last two identities therefore bound both powers for $P$ by the same product:
\[
 W(P,b)\le W(p_1,h_{p_1})^{\alpha_1}W(p_2,h_{p_2})^{\alpha_2}.
\]
It follows by H\"older's inequality that if $\Phi_Q(s;p_i)\le C$ for $i=1,2$, then $\langle W(P,b)\rangle_Q\le C$.

We now replace $b$ by the harmonic mean $h_P$. The exponents defining $W(P,b)$ are $P(1/b-1/P)$ and $-P'(1/b-1/P)$. Both $P$ and $P'$ are at least one, so $\log W(P,b)\ge|\log s|\,|1/b-1/P|$. Averaging and using Jensen's inequality gives
\[
 |\log s|\,|1/b-1/h_P|
 \le \langle\log W(P,b)\rangle_Q\le\log C.
\]
Let $K=\max\{P_+,(P')_+\}$. Replacing $b$ by $h_P$ changes the logarithm of the first power by $P(1/h_P-1/b)\log s$, whose absolute value is at most $K\log C$. The same bound holds for the conjugate power, since $1/h_P'-1/b'=-(1/h_P-1/b)$. Thus $W(P,h_P)\le C^K W(P,b)$, and averaging gives $\Phi_Q(s;P)\le C^{K+1}$.

By Theorem~\ref{main} and Lemma~\ref{centers}, we may choose a common $C>1$ for which the two integrals for $p_1(\cdot)$ and $p_2(\cdot)$, formed with $\Phi_Q$, are finite. On each cube, the set of $0<s<1/|Q|$ where $\Phi_Q(s;P)>C^{K+1}$ is contained in the union of the two sets where $\Phi_Q(s;p_i)>C$. Taking lengths and integrating proves boundedness on $L^{P(\cdot)}$.
\end{proof}

For constant shifts the comparison is sharper: one can use the arithmetic means in the statement of Theorem~\ref{main} without any change in $C$.
\begin{proposition}\label{shift}
Let $1<p_-\le p_+<\infty$ and $c>0$. Write $F_Q(s;p+c)$ and $F_Q(s;p)$ for the functions in Theorem~\ref{main} formed with $p(\cdot)+c$ and $p(\cdot)$, respectively. Then
\[
 F_Q(s;p+c)\le F_Q(s;p)\qquad(Q\subset\R^n,\ s>0).
\]
Consequently, if $M$ is bounded on $L^{p(\cdot)}$, it is bounded on $L^{p(\cdot)+c}$ for every $c>0$.
\end{proposition}
\begin{proof}
We use a comparison for normalized averages. Let $r$ take values in a compact interval in $(0,\infty)$, and let $\varphi$ be positive, continuous and nondecreasing on this interval, with $\varphi(z)/z$ nonincreasing. Then, for every differentiable convex function $H$,
\begin{equation}\label{convex-shift}
 \left\langle H\!\left(\frac{\varphi(r)}{\langle\varphi(r)\rangle_Q}\right)\right\rangle_Q
 \le\left\langle H\!\left(\frac{r}{r_Q}\right)\right\rangle_Q .
\end{equation}
To see this, choose $z_0$ so that $\varphi(z_0)/z_0=\langle\varphi(r)\rangle_Q/r_Q$. Such a point exists because the right-hand side is an average of $\varphi(r)/r$ with weight $r$. Put $b=z_0/r_Q$. Where $r\le z_0$, the two monotonicity assumptions give
\[
 \frac{r}{r_Q}\le\frac{\varphi(r)}{\langle\varphi(r)\rangle_Q}\le b ;
\]
where $r\ge z_0$, both inequalities reverse. Subtract the tangent line to $H$ at $b$. The resulting function is nonincreasing to the left of $b$ and nondecreasing to its right. Thus its value at the middle expression is at most its value at $r/r_Q$ in both cases. Averaging proves \eqref{convex-shift}, since both normalized functions have mean one and hence give the same average for the subtracted line.

Take $H(z)=s^{z-1}$, which is convex for every $s>0$. For the first average in $F_Q(s;p+c)$ use $r=p$ and $\varphi(z)=z+c$ in \eqref{convex-shift}. For the other average, observe that $(p+c)'=\varphi(p')$, where
\[
 \varphi(z)=1+\frac{z-1}{1+c(z-1)},\qquad z>1 .
\]
This function is increasing. Its ratio $\varphi(z)/z$ is nonincreasing because
\[
 1-\frac{\varphi(z)}{z}
 =\frac{z-1}{z}\,\frac{c(z-1)}{1+c(z-1)}
\]
is a product of two nonnegative increasing functions. Hence \eqref{convex-shift} also applies with $r=p'$. Neither average increases, so neither does their maximum. The set of $s$ counted by Theorem~\ref{main} for $p(\cdot)+c$ is therefore contained in the corresponding set for $p(\cdot)$ on every cube. The boundedness assertion follows.
\end{proof}

Subtracting a positive constant can destroy boundedness even if the resulting exponent stays bounded away from one \cite[Theorem~1.7]{Shifts}. For a fixed bounded function $b(\cdot)$, the set of $a$ for which $\essinf(b+a)>1$ and $M$ is bounded on $L^{b(\cdot)+a}$ is therefore an upper interval whenever it is nonempty. For a pointwise multiplier of $\mathrm{BMO}(\R^n)$ with $b_->0$, the existence of such an $a$ is the result recalled in \cite[Theorem~A]{Shifts}; every larger shift consequently gives boundedness as well.

Taking $s=1/|Q|$ in Proposition~\ref{shift} and applying Lemma~\ref{endpoint}, we see that if $p(\cdot)$ satisfies the variable
Muckenhoupt condition, then so does $p(\cdot)+c$ for every $c>0$. This answers the question raised immediately after \cite[Question~4.75]{CF}.

\subsection{Smooth bumps of different heights and widths}
Work on the line. Fix a smooth nonzero function $0\le\psi\le1$ supported in $(-1/4,1/4)$. Let $I_j$ be pairwise disjoint intervals, with centers $x_j$ and lengths $w_j>0$. 
For an example in which bumps of fixed height are added to a constant
exponent, see \cite[Theorem~5.3.4]{DHHR}. Here we add bumps of varying
heights and widths to an exponent $q(\cdot)$ that need not have a limit
at infinity.

\begin{theorem}\label{bumps}
Suppose $q(\cdot)$ is locally log-H\"older continuous, $1<q_-\le q_+<\infty$, and $M$ is bounded on $L^{q(\cdot)}(\R)$. Let $0<\delta_j\le1$ and put
\[
 p(x)=q(x)+\sum_j\delta_j\psi\!\left(\frac{x-x_j}{w_j}\right).
\]
Then $M$ is bounded on $L^{p(\cdot)}(\R)$ if and only if, for some $B>0$,
\[
 \sum_j (w_j+w_j^{-1})\,e^{-B/\delta_j}<\infty .
\]
\end{theorem}

\begin{proof}
It is useful to separate the series condition into two equivalent requirements:
\begin{equation}\label{twobump}
 \sup_j\delta_j\log^+(1/w_j)<\infty,
 \qquad \sum_j w_j e^{-B_0/\delta_j}<\infty\quad\text{for some }B_0>0 .
\end{equation}
Indeed, boundedness of the individual terms $w_j^{-1}e^{-B/\delta_j}$ gives the first requirement (use $\delta_j\le1$). Conversely, if the supremum is $H$, then $w_j^{-1}\le w_j e^{2H/\delta_j}$ when $w_j<1$, and $w_j^{-1}\le w_j$ otherwise. Increasing $B_0$ by $2H$ suffices.

Suppose \eqref{twobump} holds. The first requirement implies that the added function is locally log-H\"older continuous. To see this, for points at distance $d\le1/2$ the oscillation of the $j$th summand is at most $c\delta_j\min(1,d/w_j)$. If $w_j\le d$, multiply this bound by $\log(e+1/d)$ and use \eqref{twobump}. If $d<w_j\le1$, use
\[
 \frac d{w_j}\log(1/d)=\frac d{w_j}\log(1/w_j)+\frac d{w_j}\log(w_j/d);
\]
the second term is bounded. For $w_j>1$ the uniform Lipschitz estimate suffices. At any two points at most two summands can contribute to their difference, so the local bound is uniform. It follows that $1/p(\cdot)$ is locally log-H\"older continuous. Since $|1/p-1/q|=|p-q|/(pq)\le |p-q|$, the second requirement gives
\[
 \int_{\{p\ne q\}}\exp\!\left(-\frac{B_0}{|1/p(x)-1/q(x)|}\right)\,dx
 \le \sum_j w_j e^{-B_0/\delta_j}<\infty .
\]
The perturbation result \cite[Proposition~4.8]{CF} now transfers boundedness from $L^{q(\cdot)}$ to $L^{p(\cdot)}$.

For necessity assume boundedness for both exponents. Write
$\omega_J(r)=\av{|r-r_J|}_J$ for the mean oscillation on an interval. Whenever $J$ contains the support of the $j$th bump, is contained in $I_j$, and has $3w_j/4\le |J|\le w_j$, a fixed fraction of $J$ is outside that support. Integrating the deviation of $p-q$ from its average over this fraction gives
\begin{equation}\label{oscbump}
 \omega_J(p)+\omega_J(q)\ge\omega_J(p-q)\ge c_\psi\delta_j .
\end{equation}
Take $J=I_j$. The endpoint bounds in Lemma~\ref{endpoint}, applied to $p(\cdot)$ and $q(\cdot)$, together with \eqref{zero} applied to $(r/r_J-1)\log(1/|J|)$, show that $\omega_J(r)|\log |J||$ is uniformly bounded for $r=p,q$. By \eqref{oscbump} this proves the first part of \eqref{twobump}.

We next use the integral criterion to prove the second part of \eqref{twobump}. For $r=p,q$ write
\[
 L_C(J;r)=\big|\{0<s<1/|J|:F_J(s;r)>C\}\big|,
\]
where $F_J(s;r)$ is formed with $r(\cdot)$ and $r'(\cdot)$. Choose $C>1$ large enough that the integrals in Theorem~\ref{main} are finite for both exponents. For $r=p,q$ we claim that
\begin{equation}\label{osclength}
 L_C(J;r)\ge\min\big\{|J|^{-1},\exp(-K/\omega_J(r))\big\},
 \qquad K=2\max\{p_+,q_+\}\log(1+C).
\end{equation}
If $\omega_J(r)=0$, the exponential is understood to be zero. Otherwise, if $0<s<1$ and $F_J(s;r)\le C$, applying \eqref{zero} to $(r/r_J-1)\log s$ gives
\[
 \omega_J(r)\log(1/s)\le 2r_J\log(1+C)\le K .
\]
Hence $F_J(s;r)>C$ throughout $(0,\exp(-K/\omega_J(r)))$. Intersecting with $(0,1/|J|)$ gives \eqref{osclength}.

For $|x-x_j|<w_j/16$ and $3w_j/4<v<7w_j/8$, the interval $J=Q_v(x)$ has the properties required in \eqref{oscbump}. At least one of its two oscillations is at least $c_\psi\delta_j/2$. By \eqref{osclength}, for a common $B_0$,
\[
 L_C(J;p)+L_C(J;q)\ge\min\{w_j^{-1},e^{-B_0/\delta_j}\}.
\]
The indicated sets of centers are disjoint. Integration over each set of centers contributes a factor $w_j/8$, and integration in $dv/v$ contributes $\log(7/6)$. Finiteness of the two integrals implies
\[
 \sum_j \min\{1,w_j e^{-B_0/\delta_j}\}<\infty .
\]
Only finitely many terms can equal one, so the untruncated series converges. This proves the second part of \eqref{twobump}.
\end{proof}

\section*{Acknowledgement}
The author used ChatGPT to explore mathematical ideas, develop arguments,
and revise the exposition. The author takes full responsibility for the
results and proofs presented in this paper.

\end{document}